\documentclass[11pt]{article}
\input epsf.tex
\usepackage{amssymb,latexsym,times}
\catcode`\@=11
\renewcommand\subsection{\@startsection{subsection}{2}{\z@}%
                                     {-3.25ex\@plus -1ex \@minus -.2ex}%
                                     {-0.01 mm}
                                     {\normalfont\large\bfseries}}

\catcode`\@=11
\renewcommand\subsubsection{\@startsection{subsubsection}{2}{\z@}%
                                     {-3.25ex\@plus -1ex \@minus -.2ex}%
                                     {-0.01 mm}
                                     {\normalfont\bfseries}}

\newtheorem{example}{Example}
\newtheorem{theorem}[example]{Theorem}
\newtheorem{corollary}[example]{Corollary}
\newtheorem{proposition}[example]{Proposition}
\newtheorem{lemma}[example]{Lemma}

\def\resp{{\em resp.$\ $}}

\def\proof{\medskip\noindent {\it Proof --- \ }}
\def\finex{\hfill $\Diamond$}

\def\cqfd{\hfill $\Box$ \bigskip}
\def\adots{\mathinner{\mkern2mu\raise1pt\hbox{.}
\mkern3mu\raise4pt\hbox{.}\mkern1mu\raise7pt\hbox{.}}}
\def\<{\langle\,}
\def\>{\,\rangle}

\def\SG{\mathfrak S}

\def\l{\lambda}
\def\a{\alpha}
\def\de{\delta}
\def\b{\beta}
\def\ga{\gamma}
\def\N{{\mathbb N}}
\def\Z{{\mathbb Z}}
\def\C{{\mathbb C}}

\def\Q{{\mathbb Q}}
\def\FF{{\mathbb F}}

\def\B{{\cal B}}

\def\V{{\cal V}}

\def\SS{{\cal S}}

\def\mod{\mbox{\ mod\ }}

\def\wt{{\rm wt}}
\def\g{\mathfrak g}

\def\gl{\mathfrak{gl}}
\def\Sl{\mathfrak{sl}}

\def\slchap{\widehat{\mathfrak{sl}}}

\def\A{{\cal A}}

\def\G{{\cal G}}

\def\L{\Lambda}

\def\O{{\cal O}}

\def\<{\langle}
\def\>{\rangle}

\def\CC{{\cal C}}
\def\m{\mu}

\def\le{\leqslant}
\def\ge{\geqslant}

\def\Si{\Sigma}

\def\ul{\underline{\l}}
\def\um{\underline{\m}}

\def\irr{{\rm Irr\,}}
\def\con{{\rm Con}}
\def\sy{{\rm Sy}}
\def\ssy{{\rm SSy}}
\def\pa{{\rm Pa}}
\def\es{\emptyset}
\def\ind{{\rm Ind}}
\def\z{\zeta}

\newdimen\Squaresize \Squaresize=14pt
\newdimen\Thickness \Thickness=0.5pt
 
\def\Square#1{\hbox{\vrule width \Thickness
   \vbox to \Squaresize{\hrule height \Thickness\vss
      \hbox to \Squaresize{\hss#1\hss}
   \vss\hrule height\Thickness}
\unskip\vrule width \Thickness}
\kern-\Thickness}
 
\def\Vsquare#1{\vbox{\Square{$#1$}}\kern-\Thickness}

\def\young#1{
\vbox{\smallskip\offinterlineskip
\halign{&\Vsquare{##}\cr #1}}}
 
\title{\bf Constructible characters and canonical bases}

\author{Bernard {\sc Leclerc} and Hyohe {\sc Miyachi}}

\date{}

\begin{document}
\maketitle

\vskip 1cm

\begin{abstract}\noindent
We give closed formulas for all vectors of the canonical basis
of a level $2$ irreducible integrable representation of 
$U_v(\Sl_\infty)$.
These formulas coincide at $v=1$ with Lusztig's formulas for
the constructible characters of the Iwahori-Hecke algebras of type $B$
and $D$.
\end{abstract}

\vskip 0.6cm

\section{Introduction} \label{SECT1}

In \cite{LM} we have obtained closed formulas for certain
vectors of the canonical basis of the level~$1$ Fock space
representation of $U_v(\slchap_n)$, inspired by similar
formulas for decomposition numbers occuring in the modular
representation theory of the groups $GL_m(\FF_q)$.
(See also \cite{CT} for closely related results obtained
independently.)

In this paper we give closed formulas for all vectors of the
canonical basis of an irreducible integrable representation of 
level $2$ of $U_v(\Sl_\infty)$.
These formulas are also inspired by some known results 
in the representation theory of finite Chevalley groups, namely
by Lusztig's theory of families and constructible characters for  
Weyl groups of type $B$ and $D$, and more generally for the 
corresponding Iwahori-Hecke algebras with unequal parameters 
\cite{Lu84, Lu02}. 

Let $\L=\L_k+\L_{k+r}\ (k\in\Z, r\in\N)$ be a level~$2$ dominant
integral weight of $U_v(\Sl_\infty)$ (here the $\L_i$ are
the fundamental weights).
Let $V(\L)$ and $F(\L)$ be respectively the irreducible and 
the Fock space representation with highest weight $\L$, and let
$\Phi : V(\L) \longrightarrow F(\L)$ be the natural embedding
of $U_v(\slchap_\infty)$-modules (it is unique up to scalar
multiplication).
The Fock space $F(\L)$ is endowed by construction with a
standard basis $\SS(\L)=\{s_{(\l,\m)}\}$ canonically
labelled by all pairs $(\l,\m)$ of partitions.
Let $\B(\L)$ denote the Kashiwara-Lusztig canonical basis of 
$V(\L)$ \cite{Ka2,Lu90}.

On the other hand, let $H_m = H_m(q^r;q)$ be the Iwahori-Hecke
algebra of type $B_m$ over $\C(q^{1/2})$ with parameter $q^r$
for the special generator $T_0$ and $q$ for all other generators
$T_i\ (1\le i \le m-1)$.
Let $\irr H_m = \{\chi_{(\l,\m)} \mid |\l|+|\m|=m\}$ be the
set of irreducible characters of $H_m$ and $\con H_m$ the 
set of constructible characters \cite{Lu84,Lu02}.
Then our main result states that the specialisation at $v=1$
of $\B(\L)$ can be identified to $\bigcup_m \con H_m$, in the 
following sense: if we write for a vector $b\in\B(\L)$ of
principal degree $m$,
\begin{equation}\label{eq0.1}
\Phi(b) = \sum_{(\l,\m)} \a_{(\l,\m)}^b(v) \, s_{(\l,\m)},
\end{equation}
then
\begin{equation}\label{eq0.2}
\chi_b = \sum_{(\l,\m)} \a_{(\l,\m)}^b(1) \,\chi_{(\l,\m)} 
\end{equation}
belongs to $\con H_m$, and all constructible characters are obtained
in this way.
In particular two irreducible characters $\chi_{(\l,\m)}$ and 
$\chi_{(\l',\m')}$
lie in the same family if and only if the corresponding vectors
$s_{(\l,\m)}$ and $s_{(\l',\m')}$ have the same $U_v(\Sl_\infty)$-weight.

Similar results connect the constructible characters of type $D_m$ 
and the canonical basis of the representation $V(2\L_k)$.

The paper is organized as follows.
In Section~\ref{SECT2} we calculate explicitly the canonical basis
of $V(\L)$. 
In particular we show that the number of nonzero coefficients 
$\a_{(\l,\m)}^b(v)$
in the expansion of $b\in\B(\L)$ is always a power of 2, and
that all these coefficients are powers of $v$.
It is remarkable that the combinatorics involved in these
formulas is precisely the combinatorics of Lusztig's symbols,
which were introduced to parametrize the irreducible unipotent
characters of the Chevalley groups $G(\FF_q)$ of classical type
(see \cite{Lu84}). 
In Section~\ref{SECT3} we briefly review Lusztig's  
work on constructible characters and families for Iwahori-Hecke
algebras. 
In Section~\ref{SECT4} we compare the canonical basis of $V(\L)$ 
and the constructible characters of Iwahori-Hecke
algebras of type $B_m$ and $D_m$ calculated by Lusztig \cite{Lu79,Lu84,Lu02},
and we obtain our main result.
Section~\ref{SECT5} explains the relation between this result
and a theorem of Gyoja~\cite{Gy} comparing constructible 
characters and decomposition matrices of Iwahori-Hecke algebras.
Finally Section~\ref{SECT6} discusses a possible generalization
of our results to Ariki-Koike algebras.

It is interesting to note that Brundan \cite{Br} has obtained 
similar formulas for the canonical basis of the level 0 
module $\bigwedge^m \V^* \otimes \bigwedge^n \V$, where $\V$ denotes the
vector representation of $\Sl_\infty$ and $\V^*$ the dual representation.
His calculations were motivated by completely different problems
in the representation theory of the Lie superalgebra $\gl(m|n)$.


\section{Canonical bases} \label{SECT2}

\subsection{}
Fix $n \ge 2$ and let $\g=\Sl_{n+1}$. 
We consider the quantum enveloping algebra $U_v(\g)$ over
$\Q(v)$ with Chevalley generators 
$e_j, f_j, t_j \ (1\le j \le n)$ (see for example \cite{Ja}).
The simple roots and the fundamental weights are denoted by 
$\a_k$ and $\L_k \ (1\le k \le n)$ respectively, and the fundamental 
representations of $U_v(\g)$ by $V(\L_k)$.
Let $u_i \ (1\le i \le n+1)$ be the natural basis of the 
vector representation $V(\L_1)$, that is,
\[
e_j u_i = \de_{i,j+1} u_{i-1}, \quad
f_j u_i = \de_{ij} u_{i+1}, \quad
t_j u_i = v^{\de_{ij}-\de_{i,j+1}} u_i,
\quad (1\le i \le n+1,\ 1\le j \le n).
\]
The fundamental representation $V(\L_k)$ is obtained by $v$-deforming
the $k$th exterior power of $V(\L_1)$. 
It has a $\C(v)$-basis $\{u_\b\}$ labelled by
all sequences
\[ 
\b = (\b_1,\ldots ,\b_k),\quad 1\le \b_1<\cdots <\b_k\le n+1.
\]
For convenience, we sometimes identify such sequences with
$k$ element subsets of $\{1,\ldots , n+1\}$.
Thus we may write non ambiguously $j\in\b$, $\b \cup \{j\}$, 
and so on.  
The Chevalley generators act on $\{u_\b\}$ by
\begin{eqnarray}
e_j u_\b &=&
\left\{\matrix{0  &\mbox{if}\ j+1\not\in\b \ \mbox{or}\ j\in\b,\cr
               u_{\ga} &\mbox{otherwise, where}\ 
                \ga = (\b\setminus \{j+1\}) \cup \{j\};   \cr
}\right.
\\[2mm]
f_j u_\b &=&
\left\{\matrix{0  &\mbox{if}\  j+1\in \b \ {\rm or}\ j\not\in \b,\cr
                u_\ga &\mbox{otherwise, where}\ 
                \ga = (\b\setminus \{j\}) \cup \{j+1\};   \cr
}\right.
\\[2mm]\label{act-t}
t_j u_\b &=&
\left\{\matrix{vu_\b  &\mbox{if}\  j\in\b,\cr
               v^{-1} u_\b &\mbox{if}\ j+1\in\b,\cr
               u_\b &\mbox{otherwise}.
}\right.
\end{eqnarray}
As is well known, since $V(\L_k)$ is a minuscule representation
$\{u_\b\}$ is nothing but the canonical basis $\B(\L_k)$
of $V(\L_k)$.
The highest weight vector is $u_{\b^k}$, where
$\b^k := (1,2,\ldots ,k)$.

\subsection{}
Let now $\L = \L_k + \L_{k+r}\ (1\le k \le k+r \le n)$ be a sum 
of two fundamental weights.
Let $V(\L)$ be the irreducible $U_v(\g)$-module with highest weight
$\L$ and set $F(\L)=V(\L_{k+r})\otimes V(\L_k)$.
We have a canonical embedding of $U_v(\g)$-modules 
$\Phi : V(\L) \longrightarrow F(\L)$ defined by mapping
the highest weight vector $u_\L$ of $V(\L)$ to 
$u_{\b^{k+r}}\otimes u_{\b^k}$.

The basis 
$\SS(\L) = \{u_\b \otimes u_\ga \mid u_\b \in \B(\L_{k+r}), \ 
u_\ga \in \B(\L_k) \}$ of $F(\L)$
will be called the standard basis.
It is labelled by all symbols
\begin{equation}\label{symbol}
S=\pmatrix{\b \cr \ga} = 
\pmatrix{\b_1,\ldots ,\b_{k+r} \cr \ga_1,\ldots ,\ga_k}
\end{equation}
with $1\le \b_1 < \cdots < \b_{k+r} \le n+1$ and 
$1\le \ga_1 < \cdots < \ga_k \le n+1$. 
We shall write for short $u_S=u_\b\otimes u_\ga$.
The symbol attached to the highest weight vector is denoted by
\begin{equation}\label{symb}
S_0=\pmatrix{\b^{k+r} \cr \b^k} = 
\pmatrix{1,\ldots ,k+r \cr 1,\ldots ,k}\,.
\end{equation}
The action of the Chevalley generators on this basis is obtained
via the comultiplication of $U_v(\g)$, namely,
\begin{eqnarray}
f_j(u_\b\otimes u_\ga) &=& u_\b \otimes f_j u_\ga + f_j u_\b \otimes t_j
u_\ga,\label{fsymbol}\\
e_j(u_\b\otimes u_\ga) &=& e_j u_\b \otimes u_\ga + t_j^{-1} u_\b \otimes e_j
u_\ga,\\
t_j(u_\b\otimes u_\ga) &=& t_j u_\b \otimes t_j u_\ga. 
\end{eqnarray}

\subsection{}\label{cristal}
Let $A$ be the subsring of $\Q(v)$ consisting of all rational
functions regular at $v=0$.
Let $L$ be the $A$-lattice of $F(\L)$ spanned by $\SS(\L)$.
Since crystal bases are compatible with tensor products \cite{Ka1},
it is clear that $(\SS(\L),L)$ is a crystal basis of $F(\L)$.
Moreover, it is easy to see that the connected component of
$u_{S_0}$ in the crystal graph of $F(\L)$
is the subgraph with vertices $u_S$ where $S$ is as in (\ref{symb})
with $\b_i \le \ga_i$ for $1\le i \le k$ \cite{KN}.
Such symbols will be called standard.
Clearly, this is the same as saying that the two rows of $S$
form the two columns of a semistandard Young tableau.

\subsection{} 
Let $U^{-}_v(\g)$ be the subalgebra of $U_v(\g)$
generated by the $f_i$'s.
Let $x\mapsto\overline{x}$ denote the ring automorphism
of $U^{-}_v(\g)$ defined by
\[
\overline{f_i}=f_i,\qquad \overline{v}=v^{-1}.
\]
This induces a $\C$-linear map $u\mapsto \overline{u}$ on $V(\L)$
given by
\[
\overline{(xu_{S_0})}=\overline{x} u_{S_0},\qquad (x\in U_v^-(\g)).
\]
By \ref{cristal}, the canonical basis (or lower global basis) of
$V(\L)$ is parametrized by the set of standard symbols, 
and the element $b_S$ attached to the symbol $S$ is 
characterized by
\begin{equation}\label{canon}
\overline{b_S}=b_S
\quad \mbox{and} \quad 
\Phi(b_S) \equiv u_S \mbox{\ mod\ } vL. 
\end{equation}

\subsection{} \label{paires}
Let $S={\b\choose\ga}$ be a standard symbol.
We define an injection $\psi : \ga \longrightarrow \b$ such
that $\psi(j)\le j$ for all $j\in\ga$. 
To do so it is enough to describe the subsets
\[
\ga^l = \{j\in\ga \mid \psi(j)=j-l\},\qquad (0\le l \le n).
\]
We set $\ga^0 = \ga \cap \b$ and for $l\ge 1$ we put
\[
\ga^l = \{j\in\ga - (\ga^0 \cup \cdots \cup \ga^{l-1}) \mid 
j-l\in \b-\psi(\ga^0 \cup \cdots \cup \ga^{l-1})\}.
\]
Observe that the standardness of $S$ implies that $\psi$ is
well-defined.
\begin{example}\label{ex1}
{\rm
Take 
\[
S=\pmatrix{1 & 3 & 5 & 8 & 9 \cr
           3 & 6 & 7 & 10}\,.
\]
Then 
\[
\ga^0 = \{3\},\ \ga^1 = \{6,10\},\ \ga^2=\cdots =\ga^5=\emptyset,
\ \ga^6 = \{7\}.
\]
Hence 
\[
\psi(3)=3,\ \psi(6)=5,\ \psi(7)=1,\ \psi(10)=9.
\] 
\finex}
\end{example}

The pairs $(j,\psi(j))$ with $\psi(j)\not = j$ (that is, 
$j\not \in \b\cap\ga$) will be called the pairs of $S$.
Given a standard symbol $S$ with $p$ pairs, we denote by
$\CC(S)$ the set of all symbols obtained from $S$ by 
permuting some pairs in $S$ and reordering the rows. 
We consider $S$ itself as an element of $\CC(S)$, hence
$\CC(S)$ has cardinality $2^p$.
For $\Si\in\CC(S)$ we denote by $n(\Si)$ the number of
pairs permuted in $S$ to obtain $\Si$.
\begin{example}
{\rm We continue Example~\ref{ex1}.
There are $3$ pairs in $S$, namely $(6,5),\ (7,1),\ (10,9)$.
The $\Si\in\CC(S)$ are given below, together with the corresponding $n(\Si)$: 
\begin{center}
\begin{tabular}{|c|c|c|c|}
\hline $\Si$ & $n(\Si)$ & $\Si$ & $n(\Si)$ \\ [2mm]
\hline $\pmatrix{1 & 3 & 5 & 8 & 9 \cr3 & 6 & 7 & 10}$ & 0 &
$\pmatrix{1 & 3 & 6 & 8 & 9 \cr3 & 5 & 7 & 10}$ & 1 \\[4mm]
\hline $\pmatrix{3 & 5 & 7 & 8 & 9 \cr1 & 3 & 6 & 10}$ & 1 &
$\pmatrix{1 & 3 & 5 & 8 & 10 \cr3 & 5 & 6 & 9}$ & 1 \\[4mm]
\hline $\pmatrix{3 & 6 & 7 & 8 & 9 \cr1 & 3 & 5 & 10}$ & 2 &
$\pmatrix{1 & 3 & 6 & 8 & 10 \cr3 & 5 & 7 & 9}$ & 2 \\[4mm] 
\hline $\pmatrix{3 & 5 & 7 & 8 & 10 \cr1 & 3 & 6 & 9}$ & 2 &
$\pmatrix{3 & 6 & 7 & 8 & 10 \cr1 & 3 & 5 & 9}$ & 3 \\[4mm]
\hline     
\end{tabular}
\end{center}
\finex}
\end{example}
We can now state the main result of this section.
\begin{theorem}\label{theo}
Let $S$ be a standard symbol and let $b_S$ be the
element of the canonical basis of $V(\L)$ such that
$\Phi(b_S) \equiv u_S$ mod $vL$.
We have
\[
\Phi(b_S) = \sum_{\Si\in\CC(S)} v^{n(\Si)}\,u_\Si\,.
\]
\end{theorem}

\proof
Set $S={\b\choose\ga}$.
To simplify notation, we write $b_S$ instead of $\Phi(b_S)$
throughout this proof.
We proceed by induction on the principal degree of
$b_S$, that is, on 
\[
d=\sum_i\b_i + \sum_j\ga_j - {k+1\choose 2} - {k+r+1\choose 2}.
\]
Clearly, $d=0$ if and only if 
\[
S=S_0={\b^{k+r}\choose \b^k}
\]
is the symbol attached to the highest weight vector of $F(\L)$.
In this case we have $C(S_0)=\{S_0\}$ and $b_{S_0}=u_{S_0}$, so the
statement is true.

Otherwise, we can find $i\ge 2$ in $S$ such
that $\{i,i-1\}\cap \b = \{i\}$ or $\{i,i-1\}\cap \ga =\{i\}$.
Let $j$ be the smallest of these $i$'s.

(a) Suppose that $j \in \b\cap\ga$ and $j-1\not\in S$.
We denote by $S'$ the symbol obtained by changing in $S$ the
two occurences of $j$ into $j-1$. Clearly $S'$ is again 
standard, and the pairs of $S'$ are the same as those of $S$.
By induction we may assume that
\begin{equation}\label{induction}
b_{S'}= \sum_{\Si'\in\CC(S')} v^{n(\Si')}\,u_{\Si'}\,.
\end{equation}
The element $j-1$ occurs in both rows of each $\Si'\in\CC(S')$.
By (\ref{fsymbol}), it follows that
\[
f_{j-1}^{(2)}u_{\Si'} = u_\Si
\]
where $\Si$ is obtained from $\Si'$ by changing the two occurences
of $j-1$ into $j$.
Therefore,
\[
f_{j-1}^{(2)} b_{S'} =  \sum_{\Si'\in\CC(S')} v^{n(\Si')}\,u_{\Si}
                     =  \sum_{\Si\in\CC(S)} v^{n(\Si)}\,u_{\Si}.
\]
In particular, $f_{j-1}^{(2)} b_{S'} \equiv u_S \mod vL$.
Since $\overline{f_{j-1}^{(2)}}=f_{j-1}^{(2)}$, it follows from
(\ref{canon}) that $f_{j-1}^{(2)} b_{S'} = b_S$, and the result
is proved in this case.

(b) Suppose that $j \in \b\cap\ga$ and that $j-1$ occurs in 
one of the two rows of $S$ (it cannot occur in both rows by definition
of $j$).
We denote by $S'$ the symbol obtained by changing 
$j$ into $j-1$ in the other row. 
$S'$ is again standard.
By induction we may assume that (\ref{induction}) holds.
For $\Si'\in\CC(S')$ we have $f_{j-1} u_\Si = u_{\Si'}$
where $\Si$ is obtained by changing $j-1$ into $j$
in the row of $\Si'$ that does not contain $j$.
Indeed, if this row is the bottom row by (\ref{fsymbol})
no power of $v$ occurs in $f_{j-1} u_\Si$, and if
this is the top row then the bottom row has both $j$
and $j-1$, so the contribution of $t_{j-1}$ applied
to this row is $v^{1-1}=1$.
As in (a), it follows that $f_{j-1}b_{S'} = b_S$.
On the other hand we also have
\[
f_{j-1}b_{S'} =  \sum_{\Si'\in\CC(S')} v^{n(\Si')}\,u_{\Si}.
\]
Let us now compare the pairs of $S'$ and $S$. 
In $S'$ we have $1,\ldots ,j-1$ in both rows, and $j$ 
must be in the top row. 
Then, either $j$ does not belong to a pair of $S'$
or it belongs to a pair $(j+k,j)$.
The pairs of $S$ are the same as those of $S'$
in the first case, and in the second case they 
are the same except $(j+k,j)$ which becomes $(j+k,j-1)$.
Therefore in the first case we obviously have
 \[
\sum_{\Si'\in\CC(S')} v^{n(\Si')}\,u_{\Si}
= \sum_{\Si\in\CC(S)} v^{n(\Si)}\,u_{\Si}.
\]
This also holds in the second case since $j-1$ is changed
into $j$ in $\Si'$ in the row which does not contain $j$,
that is, in the row containing $j+k$. 
Hence we always have $j-1$ and $j+k$ lying in 
different rows in $\Si$.

(c) Suppose that there is a single occurence of $j$ in
$S$, and let $S'$ be the symbol obtained by changing
this $j$ into $j-1$. 
Again, $S'$ is standard, and
by induction we may assume that (\ref{induction}) holds. 

(c1) If $j-1$ occurs on both rows of $S'$, then it also
occurs in both rows of all $\Si'\in\CC(S')$, and no $j$
occurs in any $\Si'$.
Hence by (\ref{fsymbol}),
\[
f_{j-1}u_{\Si'} = u_{\Si_1}+v\, u_{\Si_2}
\]
where $\Si_1$ (\resp $\Si_2$) is obtained from $\Si'$
by changing $j-1$ into $j$ in the bottom (\resp top) row.
In particular, when $\Si'=S'$, $\Si_1=S$ because $\Si_2$
is not standard.
Hence we have again $f_{j-1}b_{S'}\equiv u_S \mod vL$
and therefore $f_{j-1}b_{S'}=b_S$.
On the other hand, the pairs of $S$ are all the pairs
of $S'$ plus the new pair $(j,j-1)$, thus we also have
\[
b_S= \sum_{\Si\in\CC(S)} v^{n(\Si)}\,u_{\Si}.
\]

(c2) If $j-1$ occurs in only one row of $S'$ then  
all $\Si'\in\CC(S')$ contain also one $j-1$ and no $j$.
It follows that $f_{j-1}u_{\Si'} = u_{\Si'}$ where
$\Si'$ is obtained by changing this $j-1$ into $j$.
Hence, $f_{j-1}b_{S'}=b_S$.
Finally, $j-1$ belongs necessarily to the top row
of $S'$ (otherwise, by the definition of $j$ and the
standardness of $S'$ we should have $1,2,\ldots j-1$
in both rows of $S'$ and we would be in case (c1)).
If $j-1$ does not belong to a pair of $S'$ then
the pairs of $S$ are exatly the same as those of
$S'$ and 
 \[
\sum_{\Si'\in\CC(S')} v^{n(\Si')}\,u_{\Si}
= \sum_{\Si\in\CC(S)} v^{n(\Si)}\,u_{\Si}.
\]
If $S'$ has a pair $(j+k,j-1)$, this pair becomes
$(j+k,j)$ in $S$ and all other pairs are preserved,
so the same formula still holds.
\cqfd

The proof of Theorem~\ref{theo} also shows the following
\begin{proposition}\label{monomial}
Let $\L$ be a sum of two fundamental weights.
Then each element $b$ of the canonical basis $\B(\L)$
is of the form
\[
b=f_{i_1}^{(r_1)}\cdots f_{i_s}^{(r_s)} u_\L
\]
for some $i_j\in\{1,\ldots , n\}$ and $r_j\in\{1,2\}$.
\cqfd 
\end{proposition}

Note that this Proposition can easily be proved directly by
using the fact that all $i$-strings of the crystal graph of 
$V(\L)$ have length $\le 2$ (see for example \cite{Br} 3.19). 

\subsection{} \label{pair_part} 
There is an alternative way to index the standard basis $\SS(\L)$,
namely by pairs of Young diagrams (or partitions) with shifted
content.
Let $\sy(n,k,r)$ denote the set of all symbols $S$ as in 
Eq.~(\ref{symbol}).
To $S\in\sy(n,k,r)$ we attach the pair of partitions $(\l,\mu)$ 
(written in weakly increasing order) defined by
\begin{equation}
\l_i=\b_i-i,\quad \mu_j=\ga_j-j,\qquad (1\le i\le k+r,\ 1\le j\le k).
\end{equation}
This establishes a one-to-one correspondence between
$\sy(n,k,r)$ and the set $\pa(n,k,r)$ of pairs 
$(\l,\mu)=((\l_1,\ldots,\l_{k+r}),(\mu_1,\ldots,\mu_k))$
such that 
\begin{equation}\label{part_in_rect}
0\le \l_1 \le \cdots \le \l_{k+r}\le n+1-k-r,\qquad
0\le \mu_1 \le \cdots \le \mu_k\le n+1-k.
\end{equation}
An element of $\pa(n,k,r)$ is conveniently represented by 
the pair of Young diagrams corresponding to $\l$ and $\mu$
in which the cell of $\l$ with coordinates $(i,j)$ is filled with
$i-j+k+r$ and the cell of $\mu$ with coordinates $(i,j)$ is filled with
$i-j+k$. 
Thus the symbol
\[
S=\pmatrix{1 &2 &3 &6 &8\cr 2& 3& 5} \in \sy(7,5,3)
\]
corresponds to
\[
(\l,\mu)=\left(\matrix{\young{4 & 5\cr 5 & 6 & 7\cr}}\ ,\
               \matrix{\young{1\cr 2\cr 3 & 4\cr}}\right)
\in \pa(7,5,3).
\]
Condition (\ref{part_in_rect}) is equivalent to the fact
that all cells of $\l$ and $\mu$ contain integers between
$1$ and~$n$.
\begin{example}{\rm
Take $n=2$, $k=r=1$. The correspondence is given in the table below,
where  we have denoted by $\es_i$ the empty partition regarded as the highest
weight vector in $V(\L_i)$:

\medskip
\begin{center}
\begin{tabular}{|c|c|c|c|c|c|}
\hline $S$ & $(\l,\mu)$ & $S$ & $(\l,\mu)$ & $S$ & $(\l,\mu)$ \\ [2mm]
\hline $\pmatrix{1 & 2\cr 1}$ & $(\es_2,\es_1)$ &
       $\pmatrix{1 & 3\cr 1}$ & $\left(\matrix{\young{2\cr}},\es_1\right)$ &
       $\pmatrix{1 & 2\cr 2}$ & $\left(\es_2,\matrix{\young{1\cr}}\right)$ 
\\[4mm] 
\hline $\pmatrix{2 & 3\cr 1}$ & $\left(\matrix{\young{1\cr 2\cr}},\es_1\right)$ &
       $\pmatrix{1 & 3\cr 2}$ & $\left(\matrix{\young{2\cr}},\matrix{\young{1\cr}}\right)$ &
       $\pmatrix{1 & 2\cr 3}$ & $\left(\es_2,\matrix{\young{1 &2\cr}}\right)$ 
\\[4mm] 
\hline $\pmatrix{2 & 3\cr 2}$ & $\left(\matrix{\young{1\cr 2\cr}},\matrix{\young{1\cr}}\right)$ &
       $\pmatrix{1 & 3\cr 3}$ & $\left(\matrix{\young{2\cr}},\matrix{\young{1&2\cr}}\right)$ &
       $\pmatrix{2 & 3\cr 3}$ & $\left(\matrix{\young{1\cr 2\cr}},\matrix{\young{1 &2\cr}}\right)$ 
\\[4mm] 
\hline 
\end{tabular}
\end{center}
\finex}
\end{example} 

\subsection{}\label{SSECT2.7}
Theorem~\ref{theo}, which is valid for any rank $n$, extends readily
to the quantum algebra $U_v(\Sl_\infty)$ associated to the doubly
infinite Dynkin diagram of type $A_\infty$, as we shall now explain.

In this case the simple roots $\a_i$ and the fundamental weights
$\L_i$ are indexed by $i\in\Z$.
Let $\L=\L_{k+r}+\L_k$ where $k\in\Z$ and $r\in\N$.
The standard basis of $F(\L)=V(\L_{k+r})\otimes V(\L_k)$
is labelled by the set $\sy(k,r)$ of all symbols 
$S={\b\choose \ga}$,
where 
\[
\b =(\b_i;\ i\in\Z,\ i\le k+r), \quad
\ga =(\ga_i;\ i\in\Z,\ i\le k )
\] 
are now semi-infinite increasing sequences satisfying
$\b_i = \ga_i = i$ for $i\ll 0$. 
Such sequences can be regarded as the semi-infinite wedges
arising in the classical construction of the fundamental
representations of $\Sl_\infty$ by Kac and Peterson
(see \cite{Ka}, \S14.9).

We have a trivial bijection 
$T: \sy(k,r) \longrightarrow \sy(k+1,r)$
given by 
\[
T(\b)_i=\b_{i-1}+1,\quad T(\ga)_j=\ga_{j-1}+1,\quad
(i\le k+1+r,\ j\le k+1).
\]
This comes from the diagram automorphism $k \rightarrow k+1$
of $A_\infty$, which implies that the representations
$F(\L_{k+r}+\L_k)$ and $F(\L_{k+r+1}+\L_{k+1})$ of $U_v(\Sl_\infty)$
are essentially the same. 

We also have a notion of standard symbol $S\in\sy(k,r)$, namely 
when $\b_i\le\ga_i\ (i\le k)$.
As in the finite rank case, standard symbols label the
crystal basis of the irreducible module $L(\L)$.
Like in \ref{paires}, we can define the pairs of a standard symbol, 
and since $\b_i$ and $\ga_i$ coincide for $i$ small enough, there
is a finite number of them, say $p$.
This yields the subset $\CC(S)\subset \sy(k,r)$ of cardinality 
$2^p$ obtained by permuting these pairs in $S$ in all possible ways,
and for $\Si\in\CC(S)$ the integer $n(\Si)$ of pairs in which
$\Si$ differs from $S$.

Having adapted in this way the notation, Theorem~\ref{theo}
holds without modification for the algebra $U_v(\Sl_\infty)$.

Alternatively, following \ref{pair_part}, we can replace
$\sy(k,r)$ by the set $\pa(k,r)$ of all pairs $(\l,\mu)$, where
\[
(\l_i;\ i\in \Z,\ i\le k+r),\quad (\mu_i;\ i\in \Z,\ i\le k), 
\]
are weakly increasing sequences of nonnegative integers with
finitely many nonzero elements.
Equivalently, an element of $\pa(k,r)$ can be regarded as 
a pair of Young diagrams with contents 
shifted by $k$ and $k+r$, but now without the restriction 
(\ref{part_in_rect}) on the sizes of $\l$ and $\mu$.
The correspondence $\sy(k,r)\longrightarrow \pa(k,r)$ is again
given by
\begin{equation}
\l_i=\b_i-i,\quad \mu_j=\ga_j-j,\qquad (i\le k+r,\ j\le k).
\end{equation}
When indexed by $\pa(k,r)$, the vectors of the standard basis 
$\SS(\L)$ will be denoted by $s_{(\l,\mu)}$.
The elements $(\l,\mu)\in\pa(k,r)$ corresponding to the standard symbols 
are characterized by
\begin{equation}\label{Klr}
\l_i\le\mu_i,\qquad (i\le k).
\end{equation}
\begin{example}{\rm Take $r=1$ and $k\in\Z$.
The vectors of the standard basis $\SS(\L)$ of weight
\[
\nu=\L_k+\L_{k+1}-\a_{k-1}-2\a_k-2\a_{k+1}-\a_{k+2}
\] 
are labelled by pairs $(\l,\mu)\in\pa(k,1)$ with 
$\l_i=\m_i=0$ for $i\le k-2$. 
Thus, ignoring the infinite common initial string of zeros, 
we can write for short
$(\l,\mu)=((\l_{k-1},\l_k,\l_{k+1}),(\mu_{k-1},\mu_k))$.
The vectors of the canonical basis $\B(\L)$ of weight $\nu$
are labelled by the following pairs $(\l,\mu)$:
\[
((0,0,0),(3,3)),\ \  ((0,0,1),(2,3)),\ \  ((0,0,2),(2,2)),
\ \  ((0,1,1),(1,3)),\ \  ((0,1,2),(1,2)).
\]
By Theorem~\ref{theo}, the expansion of these vectors on the 
standard basis is given by the columns of the matrix:
\[
\begin{array}{c|ccccc}
((0,0,0),(3,3))&1  &   &   &   &   \\
((0,0,1),(2,3))&v  &1  &   &   &   \\
((0,0,2),(2,2))&   &v  &1  &   &   \\
((0,1,1),(1,3))&   &v  &   &1  &   \\
((0,1,2),(1,2))&v  &v^2&v  &v  &1  \\
((1,1,1),(0,3))  &   &   &   &v  &   \\
((0,2,2),(1,1))  &v^2&   &   &   &v  \\
((1,1,2),(0,2))&   &   &   &v^2&v  \\
((1,2,2),(0,1))  &   &   &v  &   &v^2\\
((2,2,2),(0,0))    &   &   &v^2&   &   \\
\end{array} 
\]
\finex
}
\end{example}

\section{Constructible characters, left cell representations and families}  
\label{SECT3}

In this section we review following \cite{Lu84, Lu02} the definition
of constructible characters and families for Iwahori-Hecke algebras.

\subsection{}
Let $(W,S)$ be a finite Coxeter group, and let $H(W)$ be the 
corresponding Iwahori-Hecke algebra over $\C(q^{1/2})$
with parameters $q^{k(s)}$. Here $k(s)\in\N$ and $k(s)=k(s')$
when $s$ and $s'$ are conjugate in $W$.
The standard basis of $H(W)$ is denoted by $\{T_w\mid w\in W\}$.

The algebra $H(W)$ is semisimple, and its irreducible characters 
are in natural bijection with those of $\C W$ via the specialization
map $\chi\in\irr H(W) \mapsto \psi\in\irr \C W$ given by
\begin{equation}\label{spec}
\psi(w) = \chi(T_w)|_{q=1}.
\end{equation}

Let $\tau : H(W) \longrightarrow \C(q^{1/2})$ be the
symmetrizing trace defined by 
$\tau(T_w) = \de_{w,1}$.
Write
\begin{equation}
\tau = \sum_{\chi\in\irr\,H(W)} c_\chi^{-1}\, \chi\,.
\end{equation}
The Schur elements $c_\chi$ are known to be Laurent polynomials
in $q^{1/2}$. Moreover the lowest exponent of $q$ in 
$c_\chi$ is of the form $-a_\chi$ for some nonnegative integer
$a_\chi$.
Hence, using the bijection (\ref{spec}), we attach to each
$\psi \in \irr W$ an integer $a_\psi:=a_\chi$ called the $a$-invariant
of $\psi$.

\subsection{}
Let $I\subset S$ and let $(W_I,I)$ be the corresponding parabolic
subgroup. 
Let $H(W_I)$ be its Iwahori-Hecke algebra with parameters 
$q^{k(s)}\ (s\in I)$. 
By the above construction we can assign to $\xi\in\irr W_I$
an $a$-invariant $a_\xi$.
One can prove that if $\psi\in\irr W$ occurs as an irreducible
constituent of 
$\ind_{W_I}^W \xi$ then $a_\psi \ge a_\xi$. 
This suggests the definition of the truncated induction, which is
the linear map given for $\xi\in\irr W_I$ by
\[
{\mathbf{j}}_{W_I}^W(\xi)
=\sum_{
\psi \in \irr W,\,
a_\psi=a_\xi 
}
\< \chi\,,\, \ind_{W_I}^W\xi \>\, \psi\,, 
\]
where
$\< \psi\,,\, \ind_{W_I}^W\xi \>$ is the multiplicity of $\psi$
in the induced character $\ind_{W_I}^W\xi$.
Note that the $a$-invariants and the map
${\mathbf{j}}_{W_I}^W$ depend on the choice of parameters $q^{k(s)}$.

\subsection{}
Using the truncated induction, the constructible characters
of $W$ (or of $H(W)$) are defined inductively in the following way:
\begin{enumerate}
\item If $W=\{ 1 \}$, only the trivial character is constructible.
\item If $W \neq \{ 1 \}$, the set of
constructible characters of $W$ consists of all characters
of the form
\[
{\mathbf{j}}_{W_I}^{W}(\varphi) 
\quad{\mbox{or}}\quad 
{\mathrm{sgn}} \otimes {\mathbf{j}}_{W_I}^{W}(\varphi) ,
\]
where ${\mathrm{sgn}}$ is the sign character of $W$,
and $\varphi$ is a constructible character of $W_I$
for some proper subset $I$ of $S$.
\end{enumerate}

\subsection{}
Using the Kazhdan-Lusztig basis of $H(W)$, 
Lusztig has defined a partition of $W$ into subsets
called left cells, and has associated to each left
cell a representation of $H(W)$. 
In the equal parameter case, that is, when all $k(s)$ are equal, 
there is an identification theorem between 
constructible representations and left cell modules 
as follows:
\begin{theorem}[Lusztig]\label{cell=const}
Assume that $W$ is a finite Weyl group and that all $k(s)$
are equal.
Then, the constructible characters coincide with the
characters of the left cell representations of $H(W)$.
\end{theorem}

Lusztig conjectures a similar result in the unequal parameter case.

\subsection{} \label{fam}
In his study of irreducible characters of Chevalley groups $G(\FF_q)$, 
Lusztig has obtained a division of the irreducible unipotent 
characters of $G(\FF_q)$ into families.
This induces for the corresponding finite Weyl group $W$ a partition 
of $\irr W$ into certain subsets also called families 
(see \cite{Lu84}).

Later, Lusztig has introduced the constructible characters of $W$
in order to obtain a more direct way of describing the families
of $\irr W$. This goes as follows.
Consider the graph $\G_W$ with set of vertices $\irr W$
in which two irreducible characters $\chi$ and $\chi'$ are 
joined if and only if there exists a constructible character
$\psi$ such that $\chi$ and $\chi'$ both appear in the
decomposition of $\psi$.
Then the families of $\irr W$ are the connected components 
of $\G_W$.


\section{Constructible characters of type $B$ and $D$} 
\label{SECT4}

\subsection{}
We denote by $W_m$ the Weyl group of type $B_m$. 
The Coxeter generators $s_0,s_1,\ldots,s_{m-1}$ fall into two
conjugacy classes $\{s_0\}$ and $\{s_1,\ldots , s_{m-1}\}$.

Let $r\in\N$.
We are going to describe, following Lusztig \cite{Lu79,Lu84,Lu02}, 
the constructible characters of $W_m$ for the choice of parameters
$k(s_0)=r$, $k(s_1)=\cdots = k(s_{m-1})=1$.

\subsubsection{}
First, recall that the irreducible characters of $W_m$ are labelled by
the bipartitions $(\l,\mu)$ of~$m$. 
Let $k\ge m$ and $n\ge m-1+k+r$. 
Then, by \ref{pair_part}, $\pa(n,k,r)$ contains all bipartitions
of $m$. 
Hence the irreducible characters may be indexed by symbols
in $\sy(n,k,r)$.
More precisely, they are parametrized by all symbols 
$S={\b\choose\ga}$ such that
\[
\sum_i\b_i + \sum_j\ga_j - {k+1\choose 2} - {k+r+1\choose 2}=m.
\] 
Let $\sy(n,k,r,m)$ denote the subset of these symbols.
For $S\in\sy(n,k,r,m)$ we denote by $\chi_S$ the corresponding
element of $\irr W_m$. 

\subsubsection{}\label{412}
For $S={\b\choose\ga}\in\sy(n,k,r,m)$ we define
\[
Z=\{z_1<z_2<\cdots <z_M\} := (\b\cup\ga) - (\b\cap\ga),
\quad
\widetilde{Z} := (Z;\b\cap\ga),
\]
and we denote by $\pi$ the map $S\mapsto \widetilde{Z}$.
Note that $M+2|\b\cap\ga|=2k+r$, hence $M-r$ is even.

\subsubsection{}\label{413}
An involution $\iota$ of $Z$ is called $r$-admissible if
\begin{enumerate}
\item $\iota$ has $r$ fixed points;
\item if $M=r$ there is no further condition; otherwise
one requires that there exist two consecutive elements
$z, z'\in Z$ such that $\iota(z)=z'$ and that the restriction 
of $\iota$ to $Z-\{z,z'\}$ is an $r$-admissible involution
of $Z-\{z,z'\}$. 
\end{enumerate}

\subsubsection{}\label{414}
We now fix $\widetilde{Z}$ and we consider the set 
$\pi^{-1}(\widetilde{Z})$ of all symbols $S\in\sy(n,k,r,m)$
such that $\pi(S)=\widetilde{Z}$.
Let $\iota$ be an $r$-admissible involution of $Z$.
Define $\CC(\widetilde{Z},\iota)$ to be the set of all symbols
$S={\b\choose\ga}\in\pi^{-1}(\widetilde{Z})$ such that
for each non trivial orbit $O$ of $\iota$, $\beta$ 
contains one element of $O$ and $\ga$ the other.
Clearly, we have that $\CC(\widetilde{Z},\iota)$ has 
cardinality $2^{(M-r)/2}$.
The following is \cite{Lu02}, 22.24.
\begin{theorem}[Lusztig]\label{propLU}
For every $r$-admissible involution $\iota$ of $Z$,
the character
\[
\psi_\iota = \sum_{S\in\CC(\widetilde{Z},\iota)} \chi_S
\]
is constructible. Moreover, all constructible characters
of $W_m$ arise in this way. 
\end{theorem}

\subsection{} We now relate constructible characters to the 
canonical bases of Section~\ref{SECT2}.
Recall that we have defined in \ref{paires} the set of pairs 
of a standard symbol $S$.
\begin{lemma}\label{symb-inv}
{\rm (a)}\ Let $S\in\pi^{-1}(\widetilde{Z})$ be a standard
symbol.
The involution $\iota$ of $Z$ whose non trivial orbits are
the pairs of $S$ is $r$-admissible.

{\rm (b)}\ Conversely, let $\iota$ be an $r$-admissible involution
of $Z$. Let $S={\b\choose\ga}\in\pi^{-1}(\widetilde{Z})$ be the
symbol such that, if $\iota(z)=z'>z$ then $z\in\b$ and $z'\in\ga$.
Then $S$ is standard and the pairs of $S$ are the non trivial 
orbits of $\iota$. 
\end{lemma}

\proof
We shall use repeatedly the following easy remark: 

\smallskip
(R)\ if $S={\b\choose\ga}\in\pi^{-1}(\widetilde{Z})$ is a symbol 
and if $z\in\b$ and $z'\in\ga$ are equal, or are consecutive in $Z$ 
with $z<z'$, then $S$ is standard if and only if the symbol 
$S'={\b'\choose\ga'}$ with $\b'=\b-\{z\}$ and $\ga'=\ga-\{z'\}$ is standard.

\smallskip
(a)\ We argue by induction on the number $(M-r)/2$ of pairs of $S$.
If $M=r$ the claim is trivial, so assume that $M>r$.
Recall the notation of \ref{paires}.
Let $l\ge 1$ be minimal such that $\ga^l\not = \emptyset$.
By the minimality of $l$, for all $0<j<l$ and all $x\in\ga - \b$
we have $x-j\not \in \b-\ga$.
Let $z'\in\ga^l$ and set $z=\psi(z')=z'-l$. 
Then $z$ and $z'$ are consecutive in $Z$. Indeed, otherwise
there would exist $z''\in Z$ with $z<z''<z'$.
But this would contradict the minimality of $l$,
since if $z''\in \b$ we would have $z'-j\in \b-\ga$ for
$j=z'-z''<l$, and if $z''\in\ga$ we would have $z''-j\in \b-\ga$ for
$j=z''-z<l$. 

Let $S'={\b'\choose\ga'}$ with $\b'=\b-\{z\}$ and
$\ga'=\ga-\{z'\}$.
By (R), $S'$ is again standard, so by induction
we know that the involution $\iota'$ of $Z-\{z,z'\}$ whose non 
trivial orbits are the pairs of $S'$ is $r$-admissible.
Thus, using the definition of $r$-admissibility we see that
$\iota$ is $r$-admissible.

\smallskip
(b)\ If $M=r$, it follows immediately from (R) that $S$ is 
standard. Moreover $S$ has no pair and the claim is proved.
If $M>r$, by definition of $r$-admissibility, there
exists a pair $z, z'$ of consecutive elements of $Z$ such
that $\iota(z)=z'$, $z<z'$ and the restriction $\iota'$ of
$\iota$ to $Z-\{z,z'\}$ is $r$-admissible.
By induction the symbol $S'={\b'\choose\ga'}$ built from 
$\iota'$ as in the statement of the Lemma is standard and 
its pairs are the non trivial orbits of $\iota'$.
Hence, by (R) the symbol $S={\b\choose\ga}$ with 
$\b=\b'\cup \{z\}$ and $\ga=\ga'\cup \{z'\}$ is standard.
Moreover, since $z$ and $z'$ are consecutive in $Z$ they
form a pair in $S$ and the other pairs are the pairs of $S'$.
Thus, the pairs of $S$ are the non trivial orbits of $\iota$.
\cqfd 

It follows from Theorem~\ref{propLU} and Lemma~\ref{symb-inv} 
that the constructible characters of $W_m$ can be parametrized
by the set $\ssy(n,k,r,m)$ of standard symbols of $\sy(n,k,r,m)$.
Moreover, if $S\in\pi^{-1}(\widetilde{Z})$ is standard 
and if $\iota$ is the $r$-admissible involution corresponding 
to $S$ as in Lemma~\ref{symb-inv}, we have
\begin{equation}
\CC(S) = \CC(\widetilde{Z},\iota).
\end{equation}
Therefore comparing Theorem~\ref{propLU} and Theorem~\ref{theo}
we obtain immediately the main result of this section:
\begin{theorem}\label{theo2}
Let $k\ge m$ and $n\ge m-1+k+r$. 
Let $\g=\Sl_{n+1}$ and consider the canonical basis
\[
\B_m = \{b_S \mid S\in\ssy(n,k,r,m)\}
\]
of the principal degree $m$ component of the irreducible 
$U_v(\g)$-module $V(\L_{k+r}+\L_k)$.
For $b_S\in\B_m$ write as in Theorem~\ref{theo}
\[
\Phi(b_S)=\sum_{\Sigma\in\CC(S)} v^{n(\Si)}\,u_\Si\,.
\]
Then the set of constructible characters of $W_m$ for the parameters
$q^r,q,\ldots,q$ is obtained ``by specializing $v\mapsto 1$
in $\Phi(\B_m)$'', namely 
it consists of the
\[
\psi_S =  \sum_{\Sigma\in\CC(S)} \chi_\Si, \qquad (S\in \ssy(n,k,r,m))\,.
\]
\cqfd
\end{theorem}
Taking into account the remarks of \ref{SSECT2.7} one can 
reformulate Theorem~\ref{theo2} as a statement about
$U_v(\Sl_\infty)$, as we did in the introduction.
It is in fact more natural, since in this way one gets rid
of $k$ and $n$ which are irrelevant as long as they are
large enough.
This amounts to replace symbols by certain inductive limits of them,
as in \cite{Lu02} 22.7.

\subsection{} Recall from \ref{fam} the partition of $\irr W_m$ given by
the connected components of the graph $\G_{W_m}$.
Lusztig has shown (\cite{Lu02} 22.2, 23.1) that 
$\chi_S$ and $\chi_{S'}$ belong to the same class if and only
if the symbols $S,S'\in\sy(n,k,r,m)$ have the same content,
that is, the same elements with the same multiplicities.
Using Equation (\ref{act-t}), it is easy to deduce the following
alternative description:
\begin{corollary}
$\chi_S$ and $\chi_{S'}$ belong to the same component of $\G_{W_m}$
if and only if $u_S$ and $u_{S'}$ belong to the same weight space
of the $U_v(\Sl_\infty)$-module $F(\L)$.
\end{corollary} 

\subsection{}
Let $W'_m$ denote the Weyl group of type $D_m$.
The irreducible characters of $W'_m$ are indexed by unordered
bipartitions $\{\l,\mu\}$ of $m$, with the convention that
pairs of the form $\{\l,\l\}$ label two irreducible characters.
Taking $k$ and $n$ large enough we can equivalently index 
the elements of $\irr W'_m$ by the orbits on the set of symbols
$\sy(n,k,0,m)$
of the involution $\sharp$ exchanging the two rows
of a symbol, 
except that one point orbits label in fact two characters.  
We denote the elements of $\irr W'_m$ by $\chi_S=\chi_{S^\sharp}$
in the first case, and by $\chi_S^{I}, \chi_S^{II}$ in the
second case.

Note that all generators $s$ of $W'_m$ are conjugate, hence
all parameters $k(s)$ of the Iwahori-Hecke algebra $H(W'_m)$
must be equal.

\subsubsection{}
Retain the notation of \ref{412}, \ref{413}, \ref{414}, with
$r=0$.
\begin{theorem}[Lusztig]\label{theoLusztigD}
Suppose that $Z\not = \emptyset$ and let $\iota$ be a $0$-admissible
involution of $Z$. 
Then the character
\[
\psi_\iota = {1\over 2}\sum_{S\in\CC(\iota)} \chi_S
\]
is constructible. 
On the other hand, if $Z=\emptyset$, the characters $\chi_S^{I}$ 
and $\chi_S^{II}$ are both construcible.
Moreover, all constructible characters of $W'_m$ arise in
this way.
\end{theorem}

\subsubsection{} 
It follows from Theorem~\ref{theoLusztigD}, Theorem~\ref{theo} and 
Lemma~\ref{symb-inv} that: 
\begin{theorem}\label{theo2bis}
Let $k\ge m$ and $n\ge m-1+k$. 
Let $\g=\Sl_{n+1}$ and consider the canonical basis
\[
\B_m = \{b_S \mid S\in\ssy(n,k,0,m)\}
\]
of the principal degree $m$ component of the irreducible 
$U_v(\g)$-module $V(2\L_k)$.
For $b_S\in\B_m$ write as in Theorem~\ref{theo}
\[
\Phi(b_S)=\sum_{\Sigma\in\CC(S)} v^{n(\Si)}\,u_\Si\,.
\]
(Note that if $S^\sharp = S$ then $\Phi(b_S)=u_S$.)
Then the set of constructible characters of $W'_m$ 
is obtained by ``specializing $v\mapsto 1$
in $\Phi(\B_m)$'', namely 
it consists of the characters
\[
\psi_S = {1\over2} \sum_{\Sigma\in\CC(S)} \chi_\Si, 
\]
for $S^\sharp \not = S$, and of the characters 
$\chi_S^{I}$ and $\chi_S^{II}$ for $S^\sharp = S$.
\cqfd
\end{theorem}


\section{Relation with results of Gyoja} \label{SECT5}

\subsection{}
Let $R$ be an integral domain and assume that
$q$ and $Q$ are invertible elements of $R$.
We shall denote by $H_m(Q;q)_R$ the Iwahori-Hecke algebra of 
type $B_m$ over $R$ with parameters $Q$ for the special generator $T_0$
and $q$ for the remaining generators $T_i\ (1\le i\le m-1)$.

We denote by $S_R^{(\l,\mu)}$ the Specht $H_m(Q;q)_R$-module corresponding 
to a bipartition $(\l,\mu)$ of~$m$ \cite{DJ}.
Let $\{D_R^{(\l,\mu)}\}$ be a complete set of non-isomorphic simple 
$H_m(Q;q)_R$-modules parame\-tri\-zed as in \cite{Ar2} by the Kleshchev 
bipartitions $(\l,\mu)$ of $m$.
Let $P_R^{(\l,\mu)}$ be the projective indecomposable
$H_m(Q;q)_R$-module corresponding to a Kleshchev bipartition
$(\l,\mu)$, that is, the projective cover of $D_R^{(\l,\mu)}$. 

Note that when $H_m(Q;q)_R$ is semisimple, every bipartition
is a Kleshchev bipartition.
On the other hand, in the case where $Q=-q^r$ and $q$ has infinite
multiplicative order, the Kleshchev bipartitions are precisely
the elements of $\pa(k,r)$ corresponding to standard symbols 
described in Section~\ref{SECT2}, Eq.~(\ref{Klr}).

\subsection{}
Let $\A=\Z[q^{1/2},q^{-1/2}]$,
where $q^{1/2}$ is an indeterminate.
We consider the following modular system $(K,\O,\Bbbk)$.
Let ${\A_2}$ be the localization of $\A$ at the
prime ideal $2 \A$, let $\O$ be its completion, $K$ the fractional 
field of $\O$, and $\Bbbk=\FF_2(q^{1/2})$ the residue
field of $\O$.
Then, $K$ contains $\Z[q^{1/2}]$ and $q$ has clearly infinite 
order in $\O$ and $\Bbbk$.

For a finite-dimensional associative algebra ${\mathfrak{A}}$ let
$R_0({\mathfrak{A}})$ be the Grothendieck group of
the category ${\mathrm{mod}}$-${\mathfrak{A}}$.
We denote by $[M]$ the class of $M$ in $R_0({\mathfrak{A}})$.

Since $-q^r=q^r$ in $\Bbbk$,
the natural map $\A \rightarrow \Bbbk$ gives rise to two decomposition maps 
\[
{\mathbf{d}}^2_{r,+} : R_0(H_m(q^r;q)_K)\rightarrow 
R_0(H_m(q^r;q)_\Bbbk),
\quad
{\mathbf{d}}^2_{r,-} : R_0(H_m(-q^r;q)_K)\rightarrow 
R_0(H_m(q^r;q)_\Bbbk).
\]
Note that $H_m(q^r;q)_K$ is semisimple but $H_m(-q^r;q)_K$ is not 
semisimple in general.

Let $\FF$ be a field containing $\Z[q^{1/2},Q^{1/2}]$
where $Q$ is an indeterminate.
Then $H_m(Q;q)_\FF$ is semisimple.
In particular, we know that $R_0(H_m(Q;q)_\FF)$ is isomorphic to
$R_0(H_m(q^r;q)_K)$. We denote this isomorphism by ${\mathbf{i}}$.
Also, we denote by $\mathbf{d}$ the decomposition map
\[
{\mathbf{d}} : R_0(H_m(Q;q)_\FF)\rightarrow 
R_0(H_m(-q^r;q)_K). 
\]
So, we have the following commutative diagram :
\[
\matrix{
R_0(H_m(Q;q)_\FF)&\stackrel{\displaystyle\mathbf{i}}{\longrightarrow} &
R_0(H_m(q^r;q)_K)\cr\cr
\mathbf{d}\Big\downarrow & &\Big\downarrow {\mathbf{d}}^2_{r,+} \cr
R_0(H_m(-q^r;q)_K)&\stackrel{\displaystyle{\mathbf{d}}^2_{r,-}}{\longrightarrow} &
R_0(H_m(q^r;q)_\Bbbk)
}
\]

Combining Proposition~\ref{monomial} and
\cite[Corollary~3.7]{math.RT/0106185}, we deduce immediately that
\begin{proposition}\label{d=d+}
For every Kleshchev bipartition $(\l,\mu)$,
\[
{\mathbf{d}}^2_{r,-}\left([P_K^{(\l,\mu)}]\right) = [P_\Bbbk^{(\l,\mu)}].
\]
Hence, the decomposition matrices of ${\mathbf{d}}^2_{r,+}$
and ${\mathbf{d}}$ are the same and preserve their canonical
indices. \cqfd
\end{proposition}

\subsection{}
By Ariki's theorem \cite{Ar1}, the decomposition matrix of
$\mathbf{d}$ is given by the canonical basis of the 
$U_v(\Sl_\infty)$-module $V(\L_{r+k}+\L_k)$ calculated
in Section~\ref{SECT2}.

More precisely, write again $\L=\L_{r+k}+\L_k$.  
Let $\{s_{(\l,\mu)}\}$ be the standard basis of
the degree $m$ component of $F(\L)$, and for a Kleshchev bipartition
$(\l,\mu)$ let $b_{(\l,\mu)}$ be the element of the canonical basis
of $V(\L) \subset F(\L)$ such that 
$b_{(\l,\mu)} \equiv s_{(\l,\mu)} \mod vL$.
Let $\overline{F(\L)}$ and $\overline{V(\L)}$ denote the
$\Sl_\infty$-modules obtained by specializing $v$ to $1$
in $F(\L)$ and $V(\L)$ using the $\Z[v,v^{-1}]$-lattices
spanned by $\{s_{(\l,\mu)}\}$ and $\{b_{(\l,\mu)}\}$ respectively.
By abuse of notation we continue to write $\{s_{(\l,\mu)}\}$ 
and $\{b_{(\l,\mu)}\}$ for the images of these bases in
$\overline{F(\L)}$ and $\overline{V(\L)}$.

Then the complexified Grothendieck group $R_\C(H_m(Q;q)_\FF)$ is 
isomorphic to the degree $m$ component of 
$\overline{F(\L)}$, the basis $\{[S_\FF^{(\l,\mu)}]\}$ being mapped
to $\{s_{(\l,\mu)}\}$.
Similarly, $R_\C(H_m(-q^r;q)_K)$ is 
isomorphic to the degree $m$ component of
$\overline{V(\L)}$, the basis $\{[P_K^{(\l,\mu)}]\}$ being
mapped to $\{b_{(\l,\mu)}\}$. 
Finally, the map $\mathbf{d}$ corresponds to the natural 
homomorphism of $\Sl_\infty$-modules from $\overline{F(\L)}$ to 
$\overline{V(\L)}$. 

\subsection{}
In the equal parameter case, the decomposition map 
${\mathbf{d}}^2_{1,+}$ was first investigated by
Gyoja \cite{Gy} in terms of Kazhdan-Lusztig cell representations.
Taking into account Theorem~\ref{cell=const}, he proved the following 
theorem:

\begin{theorem}[Gyoja]\label{th-gyoja}
The decomposition matrix of ${\mathbf{d}}^2_{1,+}$ is equal
to the matrix whose columns give the expansion of the constructible
characters of $H_m(q;q)$ in terms of the irreducible ones.  
\end{theorem}

Thus we see that using Gyoja's theorem, Ariki's theorem and
Proposition~\ref{d=d+}, we obtain a more conceptual proof of 
Theorem~\ref{theo2} in the equal parameter case, which does not use 
explicit combinatorial calculations.

On the other hand, using our approach and Ariki's theorem, we
obtain that Theorem~\ref{th-gyoja} also holds in the unequal parameter case:
\begin{theorem}
The decomposition matrix of ${\mathbf{d}}^2_{r,+}$ is equal
to the matrix whose columns give the expansion of the constructible
characters of $H_m(q^r;q)$ in terms of the irreducible ones.  
\end{theorem}

Note that recently, Bonnaf\'e and Iancu \cite{BI} have determined
the left cells of $H_m(q^r;q)$ in the asymptotic case $r\gg m$.
They proved that in this case the characters supported by the left
cells are irreducible and coincide with the constructible 
characters.


\section{Cyclotomic algebras}  \label{SECT6}

The results of Gyoja \cite{Gy} together with some recent work
of Rouquier \cite{Ro} provide a way of generalizing the definition
of families of characters of a Weyl group to complex
reflexion groups and their cyclotomic Hecke algebras.
For the groups $W=G(d,1,m)=\Z_d\wr\SG_m$, these generalized
families have been explicitly described by Brou\'e and Kim
\cite{BK}.

\subsection{}
Let $d, m\in \N^*$ and $\mathbf{r}=(r_0,\ldots,r_{d-1})\in \Z^d$.
Let $\FF$ be a field containing $\A$ and a primitive $d$th root
of unity $\z$.
We denote by $H_{(d,m,\mathbf{r})}$ the unital associative algebra 
over $\FF$ generated by $T_0,\ldots , T_{m-1}$ subject to the relations
\[
\begin{array}{c}
(T_i+1)(T_i-q)=0 {\mbox{ \ for \ }} 1\leq i \leq m-1,\\[2mm]
\prod_{j=0}^{d-1}(T_0-\z^jq^{r_j}) =0, \\[2mm]
T_0T_1T_0T_1=T_1T_0T_1T_0,\\[2mm]
T_{i+1}T_iT_{i+1}=T_iT_{i+1}T_i
{\mbox{ \ for \ }} 1 \leq i \leq m-2,\\[2mm]
T_iT_j=T_jT_i {\mbox{ \ for \ }} 0 \leq i < j-1 \leq m-2.
\end{array}
\]
Note that the Hecke algebra $H_m(q^r;q)_\FF$ of type $B_m$
is isomorphic to $H_{(2,m,\mathbf{r})}$ with $\mathbf{r}=(r+k,k)$ 
for any $k\in\Z$.
The algebras $H_{(d,m,\mathbf{r})}$ have been introduced independently
by Ariki and Koike \cite{AK} and Brou\'e and Malle \cite{BM}. 

The algebra $H_{(d,m,\mathbf{r})}$ is semisimple and its 
irreducible modules are naturally labelled by $d$-tuples of
partitions $\ul=(\l^{0},\ldots ,\l^{(d-1)})$ with
$\sum_j |\l^{(j)}| = m$ \cite{AK}.
Denote the irreducible character of $H_{(d,m,\mathbf{r})}$ indexed
by $\ul$ by $\chi_{\ul}$.
To a $d$-partition $\ul$ one can associate a 
$d$-tuple of Young diagrams whose cells are filled with their 
content shifted by the parameters $r_j$.
More precisely, 
the cell with row number $i$ and column number 
$j$ belonging to the $s$th Young diagram is filled with the 
integer $i-j+r_s$.
Reading all the cells of $\ul$ one obtains a multiset
of integers called the content of $\ul$ that we shall denote
by $c(\ul)$.
\begin{theorem}[Brou\'e-Kim]\label{theoBK}
The characters $\chi_{\ul}$ and $\chi_{\um}$ belong to the 
same Rouquier family of $\irr H_{(d,m,\mathbf{r})}$ if and
only if $c(\ul)=c(\um)$. 
\end{theorem}

\subsection{}
To the same data we associate the $\Sl_\infty$-weight 
$\L=\sum_{j=0}^{d-1} \L_{r_j}$ and the $U_v(\Sl_\infty)$-modules
$F(\L)$ and $V(\L)$. Here, $V(\L)$ is the irreducible integrable 
module with highest weight $\L$ and 
$F(\L)=\otimes_{j=0}^{d-1}V(\L_{r_j})$.
Moreover, as before for $d=2$, the standard basis of $F(\L)$
is indexed in a natural way by all $d$-tuples of Young diagrams
with contents shifted by $r_0,\ldots , r_{d-1}$.
We shall denote this basis by $\{s_{\ul}\}$.
The $U_v(\Sl_\infty)$ weight of $s_{\ul}$ is equal to
\[
\wt s_{\ul} = \L - \sum_{j\in\Z} c_j \a_j\,
\]  
where $c_j$ denotes the number of elements equal to $j$
in $c(\ul)$ and the $\a_j$ are the simple roots of $\Sl_\infty$.
From this and Theorem~\ref{theoBK} it is easily deduced that
\begin{proposition}
The characters $\chi_{\ul}$ and $\chi_{\um}$ belong to the 
same Rouquier family of $\irr H_{(d,m,\mathbf{r})}$ if and
only if $s_{\ul}$ and $s_{\um}$ belong to the same weight
space of $F(\L)$.
\cqfd
\end{proposition}

\subsection{}
By analogy with the case $d=2$, one can then consider the
canonical basis of $V(\L)$ and its expansion on the standard
basis of $F(\L)$, 
which can be calculated via a simple algorithm \cite{LT}.
Attached to each element 
\[
b=\sum_{\ul} \a_{\ul}^b(v)\,u_{\ul}
\]
of this basis, we have a certain character of
$H_{(d,m,\mathbf{r})}$:
\[
\psi=\sum_{\ul} \a_{\ul}^b(1)\,\chi_{\ul}\,.
\]  
It would be interesting to understand whether these characters
are good analogues for $G(d,1,m)$ of the constructible or left
cell characters of the Weyl groups of type $B_m$.

\bigskip
\centerline{\bf Acknowledgements}

\bigskip
We are grateful to M. Brou\'e, M. Geck and R. Rouquier  
for very stimulating discussions.


\bigskip
\small

\noindent
\begin{tabular}{ll}
{\sc B. Leclerc} : &
Laboratoire de Math\'ematiques Nicolas Oresme,\\
& Universit\'e de Caen, Campus II,\\
& Bld Mar\'echal Juin,
BP 5186, 14032 Caen cedex, France\\
&email : {\tt leclerc@math.unicaen.fr}\\[5mm]
{\sc H. Miyachi} :&
IH{\'E}S,\\
&Le Bois-Marie, 35, route de Chartres,\\
&F-91440 
Bures-sur-Yvette, France\\
&email : {\tt miyachi@ihes.fr}
\end{tabular}

\end{document}